\documentclass{article}
\pdfoutput=1

\usepackage{amsmath,amssymb}
\usepackage{tikz}

\usepackage[margin=1in]{geometry}

\newcommand{\n}{\noindent}

\newtheorem{theorem}{Theorem}
\newtheorem{claim}{Proposition}

\newcommand{\bray}{\begin{array}}
\newcommand{\eray}{\end{array}}

\newcommand{\beq}{\begin{equation}}
\newcommand{\eeq}{\end{equation}}

\newcommand{\benum}{\begin{enumerate}}
\newcommand{\eenum}{\end{enumerate}}

\newcommand{\bc}{\begin{center}}
\newcommand{\ec}{\end{center}}

\newcommand{\bfig}{\begin{figure}[htb] \centering}
\newcommand{\efig}{\end{figure}}

\newcommand{\bitem}{\begin{itemize}}
\newcommand{\eitem}{\end{itemize}}

\newcommand{\rr}{{\cal{R}}_0}

\renewcommand{\hbar}{\bar{h}}

\newcommand{\vbar}{\bar{v}}

\newcommand{\zbar}{\bar{z}}
\newcommand{\wbar}{\bar{w}}

\newcommand{\rhobar}{\bar{\rho}}
\newcommand{\sigmabar}{\bar{\Sigma}}

\newcommand{\igrs}[1]{\includegraphics[scale=0.6]{#1}}

\title{Using asymptotics for efficient stability determination in epidemiological models}
\author{Glenn Ledder}
\date{October 29, 2023}

\begin{document}
\maketitle

\begin{abstract}
Analytical stability calculation is done to prove stability properties for systems with parameters that do not have explicit values.  For systems with three components, the usual method of finding the characteristic polynomial as $\det{J-\lambda I}$ and applying the Routh-Hurwitz conditions is reasonably efficient.  For larger systems of four to six components, the method is impractical, as the calculations become too messy.  In epidemiological models, there is often a very small parameter that appears as the ratio of a disease-based time scale to a demographic time scale; this allows efficient use of asymptotic approximation to simplify the calculations at little cost.  Here we describe the tools and an example of efficient stability analysis, followed by a set of guidelines that are generally useful in applying the method.
\end{abstract}

\section{Introduction}

Local stability analysis is an important part of the investigation of any endemic disease model.  When possible, it is ideal to identify a small set of inequalities that divide the parameter space into regions that define the set of equilibria that are locally asymptotically stable.  Often these inequalities are simple, such as when there is a unique stable endemic disease equilibrium if the basic reproduction number $\rr$ is greater than 1 and the disease-free equilibrium is stable when $\rr<1$; however, this is not always the case.  There are some models that exhibit backward bifurcation, in which there is a region in the parameter space with both a stable and an unstable endemic disease equilibrium as well as a stable disease-free equilibrium \cite{bcf, gumel, mm}.  There are also examples of models where there are no stable equilibria for certain regions in the parameter space \cite{reuf, ledder2}.

Local asymptotic stability for a given set of parameter values can be determined by numerical evaluation of the eigenvalues, but general results require analytical methods.  The direct method is to compute the characteristic polynomial as $P(\lambda)=\det (J-\lambda I)$, where $J$ is the Jacobian matrix and $I$ the identity matrix, followed by calculation of eigenvalues.  This method requires complicated calculations; hence, it is seldom used for systems with more than two components.  Systems of three components can be done without too much difficulty by finding the characteristic polynomial in the usual way and then employing the Routh-Hurwitz conditions.  In principle, the same method can be used for higher-order models, but this has rarely been done and does not appear in the principal monographs of mathematical epidemiology; e.g., \cite{bcf, mm}.

Our ability to determine stability requirements analytically depends on the efficiency of the methods we use.  Unless the Jacobian decouples, the Routh-Hurwitz conditions are more efficient than full calculation of eigenvalues.  There are also more efficient methods for calculating the characteristic polynomial.  With these tools, it is frequently feasible to do analytical stability calculation for systems with four components, as in \cite{ledder2}.  However, the additional tool of asymptotic approximation can extend the range of problems for which analytical stability analysis is possible to systems of five components, as shown below, and occasionally even six components \cite{reuf}.  This gain comes at the minor cost of a degree of error that is insignificant when compared to error caused by uncertainty of parameter values or by neglecting minor features of disease natural history.  In the development that follows, we provide detailed guidelines for the method and an illustrative example.  In so doing, we will note the specific requirements for a model to be amenable to the method.

We begin in Section \ref{timescales} with a general discussion of the characteristics of epidemiological models with two time scales.  This is followed by a description in Section \ref{tools} of the three mathematical tools that form the basis for our method.  Section \ref{two-risk} consists of an example that works through the process in some detail, culminating in a proof of stability subject to a simple, and previously unknown, restriction on one of the parameter values.  We conclude with a set of guidelines that summarizes the ideas used for the example and a discussion of the benefits of the method.

\section{Two Time Scales in Epidemiological Models}
\label{timescales}

\bfig
\igrs{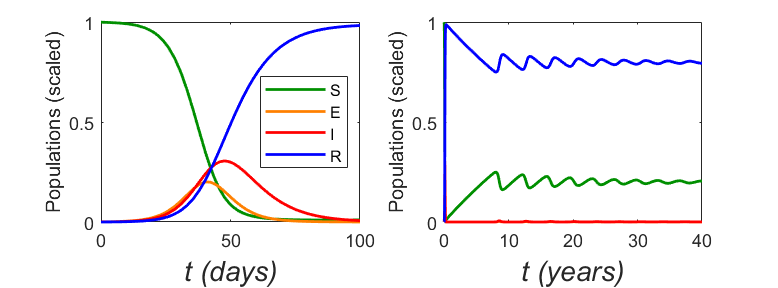}
\caption{SEIR Results on Two Time Scales}
\label{fig:SEIR}
\efig

Figure \ref{fig:SEIR} shows a typical simulation for an infectious disease model, beginning with an almost completely susceptible population.  The plots help us understand the importance of time scales in understanding epidemiological models.

The left panel shows the epidemic phase of the simulation.  We see that
\benum
\item
Infectious population fractions are significant.
\item
There is no indication of the endemic disease behavior that follows: we cannot tell whether or not there is a stable endemic disease equilibrium.
\eenum

\n In contrast, the plot on the right shows the endemic behavior.  Here,
\benum
\item
Infectious population fractions are very small.
\item
There is no indication of the epidemic behavior that preceded the endemic phase.  Everything that happened in the first 100 days occurs at the left edge of the plot.
\eenum

The difference in model behavior on the epidemic and endemic time scales is a feature of many epidemiological models that can be utilized to simplify a given task.  For investigations focused on the epidemic phase, anything that takes a long time to change the result, such as births and natural deaths, can be ignored.  For investigations focused on the endemic phase, the low populations of the infected classes can be exploited through asymptotic approximation.  This will require that the variables used to represent uninfected classes be scaled as population fractions, while the variables used to represent infected classes will need to be rescaled under the assumption that they are vanishingly small in some asymptotic limit.

\section{Mathematical Tools}
\label{tools}

The method we will employ for equilibrium point analysis is based on a combination of three mathematical tools that we describe here: an efficient scheme for calculating the characteristic polynomial, asymptotic approximation to eliminate terms with neglible impact, and the Routh array for constructing asymptotic approximations to the Routh-Hurwitz conditions. 

\paragraph{The Characteristic Polynomial Theorem}

The usual way to compute the characteristic polynomial for a matrix $J$ is to find the determinant of $\lambda I-J$;\footnote{$\lambda I-J$ is preferable to the more usual $J-\lambda I$ because its leading order coefficient is always positive.} however, there is a more efficient method based on a theorem from matrix theory (see \cite{charpoly}, for example):
\begin{theorem}
For an $n \times n$ matrix $J$, let $S$ be the set of all nonempty subsets $K$ of the integers $1, 2, \ldots, n$.  For each possible index set $K$, let $J_K$ be the determinant of the submatrix of $J$ that contains the entries in the rows and columns indicated by $K$.  Then the characteristic polynomial of $J$ is
\beq
P(\lambda)=\lambda^n + c_1 \lambda^{n-1} + c_2 \lambda^{n-2} + \cdots + c_{n-1} \lambda + c_n,
\eeq
where
\beq
c_m=(-1)^m \sum_{|K|=m} J_K, \qquad c_n=(-1)^n |J|.
\eeq
\end{theorem}
While this formula technically requires computation of all possible subdeterminants of J, this is still easier than computing the determinant of $\lambda I-J$.  The next tool will allow us to neglect some of those subdeterminants.

\paragraph{Asymptotic Approximation}

Asymptotic approximation is simply a matter of Taylor series expansion about $\epsilon \to 0$, or $\Gamma = 1/\epsilon \to \infty$ (see Chapter 1 of \cite{holmes}, for example).  Each coefficient in the characteristic polynomial has an asymptotic expansion of the form
\beq
c_m = k_m \Gamma^{p_m} + O(\Gamma^{p_m-1}), \qquad \Gamma \to \infty;
\eeq
we need only compute the $k_m$ and $p_m$ to determine the stability requirements to leading order.  For systems that have dynamic equations on two different time scales, there will be some rows in the Jacobian with a common factor of $\Gamma$ and some without.  The presence of factors of $\Gamma$ in some rows of the Jacobian means that some of the subdeterminants in the calculations of $c_1$ through $c_{n-1}$ will be dominated by others.  For example, suppose row 1 in the Jacobian has a factor of $\Gamma$ while rows 2 -- $n$ do not.  Then at least one subdeterminant of the form $J_{1j}$ will be $O(\Gamma)$, while all subdeterminants $J_{ij}$ with $i>1$ will be $O(1)$.

While one can retain terms beyond the leading order, if $\Gamma$ is sufficiently large then only terms that contribute to the leading order expansions of the characteristic polynomial coefficients need be retained to get a good approximation of stability requirements.  This generally means that only dominant terms need to be kept in all calculations; however, all Jacobian terms of $O(1)$ should be kept long enough to calculate the characteristic polynomial coefficients, even if they are dominated in their differential equation by a term of $O(\Gamma)$.  The example will show the necessity of this precaution.

\paragraph{The Routh Array}

The Routh-Hurwitz conditions use the algebraic signs of combinations of the coefficients of the characteristic polynomial to determine if all of its roots are in the left half of the complex plane.  The RH conditions for degrees 2 and 3 are well known; for higher degrees it is convenient to construct them from the Routh array, which is a two-dimensional array with $n+1$ rows and $\lceil n/2 \rceil$ columns \cite{routh}.\footnote{Other computational schemes could also be used, such as the Hurwitz matrix; however, the Routh array is particularly convenient when combined with asymptotic approximation.}  Let $R_{i,j}$ be the element in row $i$, column $j$ of the array.  The array is populated from a characteristic polynomial
\[ P(\lambda)=\lambda^n + c_1 \lambda^{n-1} + c_2 \lambda^{n-2} + \cdots + c_n \]
according to the following rules:
\benum
\item
The first row consists of the leading 1 from the characteristic polynomial, followed by $c_2$, $c_4$, and so on.
\item
The second row consists of the coefficients $c_1$, $c_3$, and so on, with a trailing 0 if $n$ is even.
\item
For $i>2$, the entries $R_{i,j}$ are given by
\[ R_{i,j} = \frac{R_{i-1,1}R_{i-2,j+1}-R_{i-2,1}R_{i-1,j+1}}{R_{i-1,1}}, \]
where formula terms requiring a column outside the dimensions of the array are taken as 0.
\eenum
As an example, consider the Routh array for a fourth degree polynomial, which takes the simple form\footnote{It is common to leave 0 entries as blanks.}
\beq
\begin{array}{ccc}
1 & c_2 & c_4 \\[.05in]
c_1 & c_3 &  \\[.05in]
\frac{q_1}{c_1} & c_4 \\[.05in]
\frac{q_2}{q_1} &  \\[.05in]
c_4 
\end{array},
\eeq
where
\beq
q_1=c_1 c_2-c_3, \qquad q_2=c_3 q_1-c_1^2 c_4.
\eeq
We'll see that expressing some of the entries as rational quantities, with a previously defined quanitity in the denominator and a newly defined quantity in the numerator, as has been done here, is a convenient modification to facilitate the next step.

The Routh theorem picks the Routh-Hurwitz stability conditions out of the Routh array.
\begin{theorem}[Routh]
The critical point with characteristic polynomial $P(\lambda)$ is locally asymptotically stable if and only if all entries in column 1 of the Routh array are positive.
\end{theorem}

\n For the fourth degree case, the conditions are
\beq
\label{rh4}
c_1>0, \qquad c_4>0, \qquad q_1>0, \qquad q_2>0.
\eeq

Asymptotic approximation sometimes reduces the formulas that arise in the Routh array to simpler ones.  For example, suppose $c_1=k_1\Gamma$, $c_2=k_2\Gamma^2$, and $c_3=k_3\Gamma^2$, all as $\Gamma \to \infty$ with $k_j=O(1)$.  Then the product $c_1c_2$ dominates $c_3$, and therefore $q_1 = k_1k_2\Gamma^3+O(\Gamma^2)$.  Since $k_1>0$ is already required, we can replace $q_1>0$ with the simpler $k_2>0$.  Simplifications such as this emerge naturally when the Routh array is constructed from the characteristic polynomial approximation
\[ P(\lambda) \sim \lambda^n + k_1 \Gamma^{p_1} \lambda^{n-1} + \cdots + k_n \Gamma^{p_n} \]
and only the leading order elements of the array are given.

\section{An SEIR Model with Two Risk Groups}
\label{two-risk}

We consider a model previously studied in \cite{gumel}. The analysis in this paper identifies a complicated parameter combination that is required for backward bifurcation to occur but does not attempt to find a simpler necessary condition.  Here we will show that backward bifurcation requires a mortality probability of greater than 75\%, an extreme condition that calls some of the model assumptions into question, and we will show that the unique endemic disease equilibrium when the mortality is not that high is always locally asymptotically stable.  

The model considers a disease scenario having two risk groups of changing composition, along with standard incidence and disease-induced mortality.  The schematic appears as Figure \ref{fig:tworisk}.  We euphemistically describe the subgroups of susceptibles as ``protected'' and ``unprotected,'' to be taken in the sense that ``protected'' people are less susceptible than ``unprotected'' people, for whatever reason.  Individuals are born/recruited into either of these categories, with some fixed proportion of each.  In the general case, they can also move between these categories through behavior changes.  By choosing parameters so that all new individuals are initially unprotected and protected individuals never switch to unprotected, we obtain the special case of a disease with a vaccine that reduces the risk of infection but does not eliminate it.  The protected are less susceptible to the disease by a factor of $\sigma<1$.  Given that the protected, unprotected, and total susceptible population fractions relative to the disease-free total population are $P$, $U$, and $S=P+U$, the total rate of infection can be written as $\beta QI/N$ where $Q=U+\sigma P=(1-\sigma) U+\sigma S$.  Note that this definition implies $U \le Q \le S$.  We can think of $Q$ as a ``susceptibility fraction'' that indicates the size of an unprotected population that would have the same total susceptibility in the absence of a protected subclass as the actual overall population.

\bfig
\begin{tikzpicture}[scale=1.25]   
\foreach \x\y\t in 
{2/2.8/$P$, 2/.8/$U$, 4.2/1.8/$E$, 6.4/1.8/$I$, 8.6/1.8/$R$}
{
\node () at (\x,\y) {\t};
\draw [rounded corners] (\x-.4,\y-.4) rectangle (\x+.4,\y+.4); 
}

\foreach \x\y\t in 
{2/.8/$U$, 4.2/1.8/$E$, 8.6/1.8/$R$}
{
\draw [->] (\x,\y-.4) to node [right, pos=0.5] {$\mu$\t} (\x,\y-.8);
}

\draw [-] (6.4,1.4) to node [right, pos=0.5] {$\mu I$} (6.4,1);

\draw [->] (6.4,1) to node [right, pos=0.5] {$\delta I$} (6.4,0.8);

\draw [->] (2,2.4) to node [below=5pt, pos=0.5] {$\mu P$} (2,2);

\foreach \x\y\t in 
{5.3/1.8/$\eta E$, 7.5/1.8/$\gamma I$}
{
\draw [->] (\x-.7,\y) to node [above, pos=0.5] {\t} (\x+.7,\y);
}

\draw [->] (2.4,2.8) to node [above=4pt, pos=0.5] {$\sigma \beta PI/N$} (3.8,2.0);

\draw [->] (2.4,0.8) to node [below=6pt, pos=0.5] {$\beta UI/N$} (3.8,1.6);

\draw [->] (0.2,2.8) to node [above, pos=0.5] {$(1-f)\mu$} (1.6,2.8);

\draw [->] (0.2,0.8) to node [above, pos=0.5] {$f\mu$} (1.6,0.8);

\draw [->, out=150, in=210] (1.6,1.2) to node [left=-2pt, pos=0.5] {$\Omega U$} (1.6,2.4);

\draw [->, out=330, in=30] (2.4,2.4) to node [right=-2pt, pos=0.5] {$\Psi P$} (2.4,1.2);

\end{tikzpicture}
\caption{Two Risk Group Model}
\label{fig:tworisk}
\efig

\subsection{Problem Definition}

From Figure \ref{fig:tworisk}, using total susceptible population $S$ in place of $P$ and adding a (redundant) equation for total population $N$, the model is
\begin{align}
\begin{split}
\tfrac{dS}{dT} &= \mu(1-S)-\beta Q \tfrac{I}{N}, \\
\tfrac{dU}{dT} &= f\mu+ \Psi (S-U)-(\Omega +\mu)U-\beta U \tfrac{I}{N}, \\
\tfrac{dE}{dT} &= -(\eta+\mu)I+\beta Q \tfrac{I}{N}, \\
\tfrac{dI}{dT} &= \eta E-(\gamma+\delta+\mu)I, \\
\tfrac{dR}{dT} &= \gamma I-\mu R, \\
\tfrac{dN}{dT} &= \mu(1-N)-\delta I,
\end{split}
\end{align}
with
\beq
N=S+E+I+R, \qquad Q=(1-\sigma) U +\sigma S, \qquad U \leq Q \leq S
\eeq
and with $T$ used for dimensional time in order to reserve $t$ for dimensionless time.
The $R$ equation will not be needed, leaving a system of five differential equations.  The constant birth rate coefficient has been assigned the same symbol ($\mu$) as that for the natural death rate, thereby automatically scaling the model so that the disease-free population total is $N=1$.  

A full stability analysis for systems with five variables can be done using the Routh-Hurwitz criteria, but the algebra will be obscure and close to unmanageable.  Instead, we will achieve approximate stability results in a relatively simple form by making use of asymptotic approximation.  In preparation for this, we need to scale the time and rescale the two infected class populations.  Scaling for dynamical systems in epidemiological models has been treated elsewhere in more detail \cite{ledder}; here we outline the best choices for this model with a brief discussion.

There are a variety of choices for the time scale.  We consider the natural demographic processes to operate on a long time scale of years, while the infection, incubation, and recovery processes operate on a short time scale of days or weeks. Thus, we have the slow rate $\mu$ and fast rates $\gamma$, $\delta$, $\eta$, and $\beta$, along with fast rates that combine multiple processes, such as the fast rate $\gamma+\delta+\mu$ that represents the rate by which individuals leave the infectious class by any means.  It is not immediately clear how to classify the status-change processes.  We tentatively take them to be slow rates and justify this choice later.  

It is important to choose one rate to represent each of the slow and fast dynamics.  The best rate to represent fast processes is $\gamma+\delta+\mu$ because the most important fast transition is the combination of transitions that remove individuals from the infectious class.  The only choice for the reference slow rate is $\mu$, which is nearly always the best choice even when there are alternatives.  We choose $1/\mu$, rather than $1/(\gamma+\delta+\mu)$ as the time scale for the dimensionless version of the problem because our interest is primarily in the long-term behavior rather than simulation.  In the nondimensionalization, we define our asymptotic parameter $\epsilon$ as the ratio of the principal slow rate to the principal fast rate and refer other rates to the appropriate principal rate according to whether the rate in question is fast or slow.  These choices lead to dimensionless parameters
\beq
\epsilon = \frac{\mu}{\gamma+\delta+\mu}, \quad b= \frac{\beta}{\gamma+\delta+\mu}, \quad m = \frac{\delta}{\gamma+\delta+\mu}, \quad \rho= \frac{\eta}{\gamma+\delta+\mu}, \quad \psi=\frac{\Psi}{\mu}, \quad \omega=\frac{\Omega}{\mu},
\eeq
and variable change
\beq
\frac{d}{dT} = \mu \frac{d}{dt}.
\eeq
For algebraic convenience, we also define
\beq
\label{params1}
\Sigma=\psi+\omega, \qquad \kappa=\sigma b, \qquad h=(1-\sigma)b, \qquad \bar{w}=w+1,
\eeq
where $w$ is any generic quantity.  The parameter $b$ is what would be the basic reproduction number for the disease if all susceptibles were unprotected, while $\sigma b<b$ would be the basic reproduction number if all susceptibles were protected.  The actual basic reproduction number will be between these values and will have to be determined as part of the analysis.

With these definitions, the scaled model is 
\begin{align}
\begin{split}
\epsilon E' &= -(\rho+\epsilon)E+bQ \tfrac{I}{N}, \\
\epsilon I' &= \rho E-I, \\
S' &= 1-S-\epsilon^{-1} bQ \tfrac{I}{N}, \\
U' &= f + \psi S-\sigmabar U -\epsilon^{-1} b U \tfrac{I}{N}, \\
N' &= 1-N-\epsilon^{-1} mI.
\label{second}
\end{split}
\end{align}
The asymptotic parameter $\epsilon$ appears on the left side of the  differential equations for $E$ and $I$ as a time scale parameter, which marks these as fast equations.  We have reordered the system to put the fast equations first; this is convenient but not necessary.  The parameter $\epsilon$ also appears on the right side of the $E$ equation; the corresponding term $-\epsilon E$ represents a slow process in a fast equation.  Usually terms of $O(\epsilon)$ that are dominated by a term of $O(1)$ can be classified as regular perturbation terms and neglected for a leading order analysis; however, we retain this term for now because slow processes in fast equations can not always be completely neglected. 

Asymptotics is based on the idea that terms of $O(\epsilon)$ are less important than terms of $O(1)$ in the limit $\epsilon \to 0$.  This only works if no quantities other than $\epsilon$ are arbitrarily small.  If a problem is properly scaled, the reduced problem obtained from $\epsilon =0$ should make sense on the long time scale.  Here it does not.  In the $N$ equation, for example, there is a term $-\epsilon^{-1}mI$, which appears to dominate the others if quantities other than $\epsilon$ are $O(1)$.  It therefore appears that the leading order approximation to this equation is $I=0$, which contradicts the expectation that the term $\epsilon^{-1}mI$ is of $O(\epsilon^{-1})$.  On closer inspection, we see that every term with $I$ appears to be more important than it actually should be.  Terms with $E$ are similarly misscaled.

The solution to this dilemma is to realize that the full population size is not the correct scale for the infected classes $E$ and $I$, as was apparent in Figure \ref{fig:SEIR}.  Everything makes sense if we assume these quantities are $O(\epsilon)$ rather than $O(1)$.  Thus, we rescale the infected class variables using the substitutions
\beq
E = \epsilon X, \qquad I = \epsilon Y.
\eeq
With these new variables, we obtain the correctly scaled model,
\begin{align}
\begin{split}
\epsilon X' &= -(\rho+\epsilon)X+bQ \tfrac{Y}{N}, \\
\epsilon Y' &= \rho X-Y, \\
S' &= 1-S-bQ \tfrac{Y}{N}, \\
U' &= f + \psi S-\sigmabar U -b U \tfrac{Y}{N}, \\
N' &= 1-N-mY.
\label{second}
\end{split}
\end{align}

At this stage, we can see why it was best to think of the transitions $\Omega U$ and $\Psi P$ as slow processes.  Had we assumed they were fast processes and defined $\tilde{\psi}=\Psi/(\gamma+\delta+\mu)$ instead of the dimensionless parameter $\psi$ that we actually used, this would have given us the term $\epsilon^{-1} \tilde{\psi} S$ instead of $\psi S$.  Since we cannot assume $S=O(\epsilon)$ without messing up the $S$ equation, we would have had to fix the problem by defining $\psi=\epsilon^{-1}\tilde{\psi}=O(1)$, taking us to where we ended up with the original choice of scale for $\Psi$.

Note that the total population equation, after rescaling, is
\[ N=S+R+\epsilon X+\epsilon Y. \]
To leading order, this is simply $R = N-S-O(\epsilon)$.  Hence, $S \le N$ will be required, but no required upper bounds for $X$ and $Y$ are apparent at this stage.

\subsection{The Disease Free Equilibrium and the Basic Reproduction Number}

The disease-free equilibrium has $X=Y=0$, $S=N=1$, $U=(\psi+f)/\sigmabar$.  Without need for a formal procedure, it is clear that the basic reproduction number, which would have been $bS$ if all susceptibles were equally susceptible, is instead 
\beq
\rr = bQ=\frac{h(\psi+f)+\kappa\sigmabar}{\sigmabar}.
\eeq
It is convenient to define a parameter $c$ to be the numerator of the fraction representing $\rr-1$; thus,
\beq
c = h(\psi+f)-(1-\kappa)\sigmabar. 
\eeq
Here we omit the routine demonstration via eigenvalues of the Jacobian that the DFE is stable if $c<0$ and unstable if $c>0$.  We should expect a unique EDE when $c>0$, but this will need to be demonstrated.  Similarly, there are usually no endemic disease equilibria when $c<0$, but this property has to be checked.

\subsection{EDE Existence and Uniqueness}

The equilibrium equations for $Y'=0$, $X'=0$, $N'=0$, and $S'=0$ yield
\[ \tfrac{bQ}{N} \sim 1, \qquad \rho X \sim Y, \qquad N+mY=1, \qquad S+Y \sim 1, \]
where the asymptotic operator $\sim$ indicates equality to leading order.
We define new quantities
\beq
\label{ysu}
y=\tfrac{Y}{N}, \quad s=\tfrac{S}{N}, \quad u=\tfrac{U}{N}, \quad p=\tfrac{P}{N}, \quad q=\tfrac{Q}{N}
\eeq
\beq
\label{zvw}
z=by, \quad v=\sigma z=\kappa y, \quad r=hy, \quad w=\Sigma +z. 
\eeq
In terms of these quantities, we have
\beq
\label{sofy}
bq \sim 1, \qquad \tfrac{1}{N}=1+my, \qquad s \sim 1-(1-m)y
\eeq
to leading order.
From the relations $q=\sigma s+(1-\sigma) u$ and $q=s+(1-\sigma) p$, we obtain additional useful formulas
\beq
\label{hu-hp}
hu \sim 1-\kappa+(1-m)\sigma z, \qquad hp \sim bs-1. 
\eeq

To complete the specification of the endemic disease equilibrium, we substitute the formulas for $1/N$, $s$, and $u$ into the equation $U'/N=0$, or
\beq
\label{Ueqn}
(\sigmabar+z)u=\psi s+\tfrac{f}{N},
\eeq
and simplify to obtain the quadratic equation
\beq
\label{edefn}
G(z) \equiv (1-m) \sigma z^2 + [(1-\kappa)+(1-m)\sigma \sigmabar +(1-m)(1-\sigma)\psi-(1-\sigma) mf]z +(1-\kappa)\sigmabar-h(\psi+f)=0(\epsilon),
\eeq
or
\beq
\label{edeeq}
(1-m) \sigma z^2 + [(1-\kappa)+(1-m)\sigma \sigmabar +(1-m)(1-\sigma)\psi-(1-\sigma) mf]z = c+O(\epsilon).
\eeq
Gumel showed that this model has a backward bifurcation and gave a complicated formula, involving all the model parameters, for the conditions required for it to occur \cite{gumel}.  Here we show that backward bifurcation requires unrealistically extreme parameter values.\footnote{This is not a pejorative statement, nor is it a value judgment, but is based on the assumptions used to build the model. In the discussion, we explain why the model we are using would be inappropriate with such a large disease-related death rate.}
\begin{claim}
In the asymptotic limit $\epsilon \to 0$, endemic disease equilibria for the model (\ref{second}) with $\rr \le 1$ can occur only with the combination $m>\kappa$ and $m>0.75$.
\end{claim}

To establish Proposition 1, we first note that solutions of the equation for endemic disease equilibria depend continuously on the model parameters, with the consequence that positive solutions with $c<0$ (ie., $\rr<1$) occur if and only if the second solution for the case $c=0$ (the other being $z=0$) is positive.  As the leading coefficient is positive, the non-zero solution will be positive if the middle coefficient is negative; hence, we assume that the middle coefficient of $G(z)$ is negative while $c=0$ and identify the implications of this assumption.  Multiplying through by $b$ and using $\kappa=\sigma b$, this means that we require
\[ (b-\kappa)mf>(1-m)\kappa \sigmabar+(1-m)(b-\kappa)\psi+b(1-\kappa) \]
while
\[ (1-\kappa) \sigmabar = (b-\kappa)(\psi+f). \]

Multiplying the inequality by $1-\kappa$ and substituting from the equality yields the inequality
\[ (b-\kappa)m(1-\kappa)f>(1-m)\kappa (b-\kappa)(\psi+f)+(1-m)(b-\kappa)(1-\kappa)\psi+b(1-\kappa)^2, \]
which we can rearrange to get
\[ b (1-\kappa)^2+(1-m)(b-\kappa)\psi<(m-\kappa)(b-\kappa) f .\]
We see that the inequality requires $m>\kappa$.  With the restrictions $\psi \ge 0$ and $f \le 1$ that are part of the model assumptions, the inequality reduces to
\[ b (1-\kappa)^2<(m-\kappa)(b-\kappa). \]
Rearranging yields
\[ \kappa(m-\kappa)<b(m-\kappa)-b(1-\kappa)^2=b(m-1+\kappa-\kappa^2). \]
Thus, we must have
\[ m>1-\kappa+\kappa^2=\tfrac34+ \left( \tfrac12-\kappa \right)^2 \ge \tfrac34. \]

To complete the description of existence and uniqueness properties, we consider the possibilities for solutions of (\ref{edeeq}) when $\rr>1$.
\begin{claim}
In the asymptotic limit $\epsilon \to 0$, there is a unique endemic disease equilibrium that satisfies all restrictions on the signs of the state variables whenever $\rr>1$.
\end{claim}

It is immediately clear that the equation has a unique positive solution when $c>0$, but it still needs to be shown that the solution is always biologically meaningful.  The critical requirement is $p \ge 0$, which through (\ref{hu-hp}) and (\ref{sofy}) corresponds to
\[ z^* \le \tfrac{b-1}{1-m} \equiv \hat{z}, \]
where $z^*$ is the positive root of the function $G$ defined in (\ref{edefn}).  Given that $G$ is a quadratic function that is negative at $z=0$, $0$ at $z=z^*$, and becoming infinite as $z \to \infty$, the requirement $z^* \le \hat{z}$ is equivalent to
\[ G(\hat{z}) \ge 0. \]
Replacing products $(1-m)\hat{z}$ by $b-1$, we have
\begin{align*}
G(\hat{z})= \sigma (b-1) \hat{z}+(1-\sigma b) \hat{z} &+ \sigma (b-1) \sigmabar+(1-\sigma)(b-1)\psi-(1-\sigma)mf\hat{z} \\
&+(1-\sigma b)\sigmabar-(1-\sigma)b\psi-(1-\sigma)bf,
\end{align*}
which simplifies to
\begin{align*} 
G(\hat{z}) &= (1-\sigma) \hat{z} + (1-\sigma) \sigmabar-(1-\sigma)\psi-(1-\sigma)f(b+m\hat{z}) \\
&=(1-\sigma) [\hat{z} + 1+\omega-f(b+m\hat{z}]).
\end{align*}
Given the parameter restrictions $\omega \ge 0$ and $f \le 1$, we have
\[ \tfrac{1}{(1-\sigma)}G(\hat{z}) \ge (1-m)\hat{z}+1-b=0. \]

\subsection{EDE Stability}

We begin by identifying a set of three stability requirements, given as inequalities in terms of the EDE solution $y$ and the other parameters.

\begin{claim}
An endemic disease equilibrium is locally asymptotically stable in the asymptotic limit $\epsilon \to 0$ if and only if $A,B,C>0$, where
\beq
\label{eqAB}
A=(\kappa-m)+bhu, \qquad B = A +(\kappa-m)w+mhuz +h \psi.
\eeq
\beq
\label{eqC}
C = \rhobar A+\rhobar^2(A\vbar+A\wbar-B)-\rho y A^2.
\eeq
\end{claim}

The Jacobian, given in terms of the additional symbol definitions of (\ref{params1}), (\ref{ysu}), and (\ref{zvw}), takes the compact form
\beq
J=\left( \begin{array}{crccc}
-(\rho \Gamma+1) & \Gamma & v\Gamma & r \Gamma & -y\Gamma \;\; \\
\rho \Gamma & -\Gamma & 0 & 0 & 0 \\
0 & -1 & -\vbar \;\; & -r \;\; & y \\
0 & \;\;\; -bu & \psi & -\wbar \;\; & uz \\
0 & \;\;-m & 0 & 0 & -1 \;\;
\end{array} \right),
\eeq
where $\Gamma=1/\epsilon$.  The leading order characteristic polynomial is obtained from this matrix by a combination of the characteristic polynomial theorem (Theorem 1) and asymptotic approximation.  

The coefficient $c_1$ is the negative of the trace.  Some terms are $O(\Gamma)$, while others are $O(1)$.  To leading order, we need only consider the former; hence,
\beq
c_1=\rho \Gamma+\Gamma+O(1)=\rhobar \Gamma +O(1).
\eeq

At first glance, $c_{12}$ appears to be dominated by an $O(\Gamma^2)$ term coming from $J_{12}$; however, the terms of that order cancel, giving $J_{12}=\Gamma$.  This means that the other terms of that order, coming from subdeterminants with indices 1 or 2 paired with 3, 4, or 5, enter in at leading order.  The result is
\beq
c_2 = [1+\rhobar(1+\vbar+\wbar)]\Gamma+O(1).
\eeq
Note that we needed the second term in the 11 entry even though it is subdominant to the first term.  That second term will not be needed for any other characteristic polynomial coefficients.

The coefficient $c_3$ is dominated by terms of $O(\Gamma^2)$ that come from the three subdeterminants of form $J_{12j}$.  These conform to a simple pattern:
\[ c_3 \sim -(J_{123}+J_{124}+J_{125})=\rho \Gamma^2 [g(\kappa y,1,\vbar)+g(hy,bu,\wbar)+g(-y,m,1)], \]
where
\[ g(D,E,F)=- \left| \bray{rrr} -1&1&D\\ 1&-1&0\\ 0&-E&-F \eray \right| = DE. \]
Thus,
\beq
c_3=\rho y A \Gamma^2 +O(\Gamma),
\eeq
where $A$ was defined in (\ref{eqAB}).

The full determinant calculation yields
\beq
c_5=-|J|=\rho y B \Gamma^2,
\eeq
with $B$ given in (\ref{eqAB}).

The coefficient $c_4$ is the sum of the $4 \times 4$ subdeterminants of $J$.  The three (of the total of five) that include both rows 1 and 2 contribute terms of $O(\Gamma^2)$.  Each of the three required subdeterminants is described by the template
\[ J_{12ij}=\left( \begin{array}{crcc}
-\rho \Gamma \;\; & \Gamma & \hat{A}y \Gamma & \hat{B}y\Gamma \\
\rho \Gamma & -\Gamma & 0 & 0 \\
0 & -\hat{C} & -\hat{E} \;\; & \hat{G} \\
0 & -\hat{D} & \hat{H} & -\hat{F} \;\;
\end{array} \right). \]
Cofactor expansion on the first column yields
\begin{align*} 
J_{12ij} &= -\rho \Gamma 
\left( \begin{array}{rcc}
-\Gamma & 0 & 0 \\
-\hat{C} & -\hat{E} \;\; & \hat{G} \\
-\hat{D} & \hat{H} & -\hat{F} \;\;
\end{array} \right)
-\rho \Gamma
\left( \begin{array}{rcc}
\Gamma & \hat{A} y\Gamma & \hat{B} y\Gamma \\
-\hat{C} & -\hat{E} \;\; & \hat{G} \\
-\hat{D} & \hat{H} & -\hat{F} \;\;
\end{array} \right) \\
&=-\rho y \Gamma^2
\left( \begin{array}{rrr}
0 & \hat{A} & \hat{B} \\
-\hat{C} & -\hat{E} & \hat{G} \\
-\hat{D} & \hat{H} & -\hat{F} \end{array} \right) 
=\rho y \Gamma^2 (\hat{D}\hat{A}\hat{G}+\hat{D}\hat{B}\hat{E}+\hat{C}\hat{A}\hat{F}+\hat{C}\hat{B}\hat{H}).
\end{align*}
Careful accounting of the terms in this formula from each of the three relevant submatrices yields the convenient result
\beq
c_4=c_3+c_5.
\eeq

To summarize, we have found the leading order characteristic polynomial
\beq
\label{charpoly2}
P(\lambda)=\lambda^5 +k_1 \Gamma \lambda^4 + k_2 \Gamma \lambda^3 + k_3 \Gamma^2 \lambda^2 + k_4 \Gamma^2 \lambda+k_5 \Gamma^2,
\eeq
where
\beq
k_1=\rhobar, \quad k_2=1+k_1+\rhobar(\vbar+\wbar), \quad k_3=\rho yA, \quad k_5 =\rho yB, \quad k_4=k_3+k_5,
\eeq
Note that asymptotic simplification only applies to terms within each individual coefficient.  We cannot compare the magnitudes of the different terms because there are eigenvalues of different orders.  As an example, suppose we assume $\lambda=O(1)$.  The leading order characteristic polynomial is then quadratic and has two roots.  Hence, only two of the eigenvalues are $O(1)$.  If instead we assume $\lambda=O(\Gamma)$, then the dominant terms are the first and second, leading to one eigenvalue, $-k_1 \Gamma$.  The remaining two eigenvalues are of an intermediate order.  We need not be concerned with details, as we do not need to find the eigenvalues to determine stability.

Neglecting all lower-order terms, the Routh array for the endemic disease equilibrium of the two-risk-class model, with characteristic polynomial given by (\ref{charpoly2}), takes the simplified form
\beq
\begin{array}{ccc}
1 & k_2 \Gamma & k_4 \Gamma^2 \\[.05in]
k_1 \Gamma & k_3 \Gamma^2 & k_5 \Gamma^2 \\[.05in]
\frac{q_1}{k_1} \Gamma & k_4 \Gamma^2 \\[.05in]
\frac{q_2}{q_1}\Gamma^2 & k_5 \Gamma^2 \\[.05in]
k_4 \Gamma^2 \\[.05in]
k_5 \Gamma^2
\end{array},
\eeq
where
\beq
q_1=k_1 k_2-k_3, \qquad q_2=k_3 q_1-k_1^2 k_4.
\eeq

The quantity $k_1$ is already known to be positive.  Given the requirements $q_1>0$ and $k_4>0$, the requirement $q_2>0$ shows that $k_3>0$ is also required.  Since the formula for $k_3$ is simpler than that for $q_1$, we need only check $k_3>0$ and $q_2>0$ to confirm $q_1>0$, and of course we have to check $k_5>0$.  The requirement $k_4>0$ is automatically satisfied when $k_3>0$ and $k_5>0$.  Using the simplifying notation introduced earlier, we are left with three conditions: $A>0$, $B>0$, and $C>0$, where
\[ \rho y C =q_2=k_3(k_1k_2-k_3)-k_1^2k_4=k_1(k_2-k_1)k_3-k_1^2k_5-k_3^2. \]
Substituting in the formulas for the $k$'s yields the formula given in (\ref{eqC}).

The three stability criteria given in Proposition 3 can be used to determine stability of the equilibrium solution for any specific set of parameters.  However, with careful algebra we can establish a strong result:
\begin{claim}
If $m \le \kappa$ or $m \le 0.75$, any endemic disease equilibrium for the two-risk-class model is locally asymptotically stable in the asymptotic limit $\epsilon \to 0$.
\end{claim}

The first step in establishing Proposition 4 is to obtain a tight lower bound for $C$ in terms of a formula with form similar to that of $A$ and $B$, that is, $C>(\kappa-m)C_a+C_b$, where $C_a$ and $C_b$ are nonnegative.  
We begin by multiplying (\ref{Ueqn}) by $b$, rewriting in terms of $w$, and using (\ref{hu-hp}) and $p+u=s$ to get
\[ bu\wbar-\psi=b \left( \psi s+\tfrac{f}{N} \right)-\psi=\psi hp+\tfrac{bf}{N} =h(\psi s-\psi u)+\tfrac{bf}{N}. \]
Rewriting (\ref{Ueqn}) as
\[ (\omega+\zbar)u=\psi s-\psi u+\tfrac{f}{N} \]
allows us to complete the calculation as
\[ bu\wbar-\psi=hu(\zbar+\omega)+\tfrac{b \sigma f}{N}, \]
leading to a computational result that will be needed shortly:
\beq
\label{id0}
bu\wbar-\psi \ge hu \zbar.
\eeq

From
\[ \begin{array}{rrrr}
A \vbar = \quad & (\kappa-m) \vbar & +\kappa huz & +bhu \; \\
A \wbar = A & +(\kappa-m) w && +bhuw \\
B = A & +(\kappa-m) w & +mhuz & +h\psi \;\;
\end{array}, \]
we obtain
\beq
\label{iq1}
A\vbar+A\wbar-B=(\kappa-m)(\vbar+huz)+h(bu\wbar-\psi).
\eeq
Applying (\ref{id0}) gives the result
\beq
\label{id1}
A\vbar+A\wbar-B \ge (\kappa-m)(\vbar+huz)+h^2u\zbar
\eeq
Next, we write
\begin{align*}
yA^2 &= (\kappa-m)^2y + 2(\kappa-m)bhuy+ b^2h^2u^2y
\\ &= (\kappa-m)(\kappa y-my+2huz)+bh^2u^2z. 
\end{align*}
Then $u \le q$ and $bq=1$ combine to give $bu \le 1$, so
\beq
\label{id2}
yA^2 \le (\kappa-m)(v+2huz-my)+h^2uz,
\eeq
Finally, we substitute (\ref{id1}) and (\ref{id2}) into (\ref{eqC}) to obtain the desired lower bound
\beq
\label{ieqC}
C \ge (\kappa-m) \left[ \rhobar^2 +(\rhobar^2-2\rho)huz+(\rhobar^2-\rho)v+\rho my \right]+(\rhobar^2-\rho) h^2 uz+\rhobar^2 h^2 u +\rhobar A.
\eeq

From (\ref{eqAB}) and (\ref{ieqC}), we see immediately that stability is guaranteed if $m \le \kappa$.  It remains to be shown that $m \le 0.75$ is sufficient to guarantee stability when $m>\kappa$.  The first step in this demonstration is to rewrite (\ref{ieqC}) as
\[ C \ge A+\rhobar^2 C_1+(\rhobar^2-\rho)zC_2+\rho C_3, \]
where
\beq
\label{C12}
C_1=h^2u-(m-\kappa), \qquad C_2=h^2u-(m-\kappa)(hu+\sigma), 
\eeq
\beq
\label{C3}
C_3=bhu-(m-\kappa)(1+my-huz).
\eeq
Given $h<b$, $hu+\sigma=(1-\sigma)bu+\sigma \le (1-\sigma)+\sigma=1$, and $m >\kappa$, we have $A \ge C_1$ and $C_2 \ge C_1$, leaving $C_1>0$, $C_3>0$, and $B>0$ to be confirmed.

From the formula for $C_1$ and $hu \sim 1-\kappa+(1-m)\sigma z \ge 1-\kappa$, we have
\[ C_1 \ge (b-\kappa)(1-\kappa)+\kappa-m=b(1-\kappa)+\kappa^2-m > \tfrac14-\kappa+\kappa^2=\left( \tfrac12-\kappa \right)^2 \ge 0. \]
For $C_3$, we first note that $s \ge 0$ requires $y \le 1/(1-m) \le 4$, and we also have $bhu>hu \ge 1-\kappa$.
If $m>bhu$, then
\[ 1+y(m-bhu)<1+4 \left( \tfrac34-1+\kappa \right)=4\kappa \]
and 
\[ C_3>(1-\kappa)- \left( \tfrac34-\kappa \right)(4\kappa)=1-4\kappa-\kappa^2=(1-2\kappa)^2 \ge 0. \]
If $m \le bhu$, then $1+y(m-bhu) \le 1$ and
\[ C_3>bhu-(m-\kappa)=1-m>0. \]
To show $B \ge 0$, we first note the result
\[ mhu-m+\kappa \ge m(1-\kappa)-m+\kappa=(1-m) \kappa \ge 0. \]
This allows us to eliminate the $z$ terms from (\ref{eqAB}), leaving 
\[ B \ge b(1-\kappa)+h\psi-(m-\kappa) \sigmabar, \]
Multiplying by $1-\kappa$ and using $b>h$, this becomes
\[ (1-\kappa)B > h(1-\kappa)^2+(1-\kappa)h\psi-(m-\kappa)(1-\kappa) \sigmabar. \]
With $m \le 0.75$, there is no backward bifurcation; hence, solutions occur only if $c>0$; along with the modeling restriction $f \le 1$, we then have
\[ (1-\kappa)\sigmabar \le h\psi+h. \]
Substituting into the previous inequality yields
\begin{align*}
(1-\kappa)B
& > h(1-\kappa)^2+(1-\kappa)h\psi-(m-\kappa)h\psi-(m-\kappa) h
\\ & \ge h \left[ (1-\kappa)^2- \left( \tfrac34 -\kappa \right) \right]
\\ & =h \left( \tfrac12-\kappa \right)^2 \ge 0.
\end{align*}
This completes the proof of Proposition 4.

\section{Guidelines for the Asymptotic Approach}
\label{guidelines}

The theorems for calculation of the characteristic polynomial and the Routh-Hurwitz conditions can be used on any problem.  Obtaining asymptotic approximations to the stability conditions for an equilibrium requires that some guidelines be followed so that the method can be properly employed and so that algebraic complexity is kept to a minimum.  We present a set of guidelines here; the starting point is the dimensional dynamical system model.

\benum
\item
Scale populations so that the disease-free equilibrium population total is 1.  This can often be done during the initial model construction by choosing a birth rate coefficient to match the natural death rate.  
\item
Identify one principal fast time scale and one principal slow time scale.\footnote{In mathematical epidemiology, fast processes have time scales on the order of weeks, while slow processes have time scales on the order of years, or perhaps months.  Usually the principal fast time scale is the amount of time spent in the principal infectious state and the principal slow time scale is the lifespan of the disease host.}  Use the principal slow time scale to scale the time.\footnote{The slow time scale is indicated when the primary interest is in long-term behavior.}
\item
Define $\epsilon$ as the ratio of the principal slow rate to the principal fast rate.  Make all other parameters dimensionless by reference to the appropriate principal scale; that is, rate parameters in fast processes, such as infection, should be referred to the principal fast rate, while rate parameters in slow processes, such as loss of immunity, should be referred to the principal slow rate.  When in doubt, use the reference value that makes the corresponding terms important but not dominant.
\item
Rescale ``small’’ populations by a factor of $\epsilon$.  Small populations are generally those for which the dynamics of increase and decrease are dominated by fast processes; for example, in the SEIR model, the dynamics of the latent (``exposed’’) class is due to the fast processes of infection and incubation, and the dynamics of the infectious class is due primarily to the fast processes of incubation and recovery.\footnote{The infectious class in the SI model is not small in this sense because there is no recovery process to drain away the population.}  A population can be seen to be small if it appears only in terms with factors of $\epsilon^{-1}$ relative to the other terms in the equation.  These terms look large, but only because the corresponding state variable has an inherent, albeit invisible, factor of $\epsilon$ when viewed on the long time scale.  The rescaling serves to make all variables $O(1)$, so that the importance of terms is determined entirely by the factors of $\epsilon$. If no populations meet this criterion, then asymptotics will have limited value, if any.
\item
Define auxiliary quantities to simplify formulas.  Examples include quantities that represent weighted averages of infectious or susceptible populations or quantities that combine together in one entry of the Jacobian matrix.  Because of the scaling, there will often be terms that include a variable plus 1; in these cases, a more compact notation is facilitated by the convention of using a bar over the top of any quantity to indicate that quantity plus 1.  
\item
When equilibria cannot be determined by linear equations alone, use the equilibrium relations to try to write all but one of the dependent variables in terms of a principal dependent variable $z$.  Derive the resulting polynomial equation $G(z)=0$ for that dependent variable, but do not solve it.  
Parameter range requirements for solution existence can be determined in terms of the function $G$ evaluated at critical values when direct comparison of the critical value with the equilibrium values is difficult.
\item
While not necessary, it is helpful to choose an equation order for the Jacobian that groups the fast equations at the beginning of the list.  Fast equations are those that correspond to fast (rescaled) variables.  The rescaling puts a factor of $\epsilon$ on the left side of the equation (or a factor of $\Gamma=1/\epsilon$ on the right side).
\item
Look for opportunities to simplify formulas by replacing dependent variables with versions that include an extra factor, such as using $y=\rr Y$ to replace the rescaled infectious population $Y$.  Introduction of modified dependent variables by incorporating common factors can be done before or after the determination of equilibrium relations but should come after the calculation of the Jacobian. 
\item
During the stability analysis, it is important to keep checking whether combinations of auxiliary quantities can be simplified or whether it would be helpful to introduce additional combinations.
\item
Unless the eigenvalues are immediately apparent from the Jacobian, use the appropriate Routh-Hurwitz conditions to determine stability.  In doing so, use Theorem 1 to compute the coefficients of the characteristic polynomial to leading order.  For anything over 4 variables, use the Routh array (Theorem 2) to determine the Routh-Hurwitz conditions to leading order.
\eenum

\section{Discussion}

We have now seen an example of asymptotic stability analysis facilitated by asymptotic approximation in conjunction with the simplest formula for calculating the characteristic polynomial of a matrix and a mechanism for constructing Routh-Hurwitz conditions for polynomials of any degree.  Using the characteristic polynomial coefficient theorem and the Routh array come at no cost in accuracy, while the use of asymptotic approximation comes with the cost of using an approximate result in place of an exact result.  We now discuss in more detail the advantages and disadvantages of approximation.

The main advantage of asymptotic approximation is that it can simplify complex formulas to the point where intractable ones become tractable, as in the stability demonstation for the two-risk-class model.  While this demonstration required some messy algebra, it is instructive to compare it to the alternative computation performed without benefit of asymptotics.  Using the generic form for the characteristic polynomial that includes the $c$ coefficients rather than products of the $k$'s with factors of $\Gamma$, the entry in row 5, column 1 of the Routh array is given by $q_4>0$, where $q_4$ is calculated from the successive formulas
\beq
q_1=c_1c_2-c_3, \quad q_2=c_3 q_1-c_1^2q_3, \quad q_3=c_1c_4-c_5, \quad q_4=q_2q_3-c_5q_1^2.
\eeq
If these are combined into a single formula with all parentheses removed, the result is
\beq 
\label{mess}
q_4=c_1 \left[ c_1c_2c_3c_4+c_2c_3c_5+2c_1c_4c_5-c_3^2c_4-c_1^2c_4^2-c_1c_2^2c_5-c_5^2 \right].
\eeq
This is sufficiently complicated to deter anyone from attempting analytical stability analysis for systems with five components.

Now suppose we apply asymptotics from this starting point.
When we replace the $c$'s with factors of $k$'s and $\Gamma$, we discover that the first, fourth, and fifth terms in $q_4$ are $O(\Gamma^6)$, while the other three are $O(\Gamma^5)$.  Removing the lower order terms reduces the formula to $q_4=c_1c_4q_2$.  Other Routh-Hurwitz conditions have already required $c_1$ and $q_2$ to be positive, so the condition reduces further to the single condition $c_4>0$, corresponding to the condition $k_4>0$ that appeared in row 5, column 1 of the Routh array with asymptotics.  In our development, we showed that this condition is not even binding, a result that would have been impossible to obtain from (\ref{mess}).

There is no question that the asymptotic results in the limit $\epsilon \to 0$ found here could not have been found in the more general case where $\epsilon$ is not taken to be arbitrarily small.  There does remain the question of whether results that require an asymptotic limit have general value.  In problems such as the two-risk class model, where stability is guaranteed for any reasonable parameter values ($m \le 0.75$), asymptotic stability as $\epsilon \to 0$ must at least hold for $\epsilon$ values below some finite threshold.   In epidemiological models where the slow and fast time scales are as different as the duration of a typical human disease and the typical human life span, actual values of $\epsilon$ are quite small, so there is a very high likelihood that the results hold for any realistic values of $\epsilon$.

We stated earlier that the requirement $m \le 0.75$ is not really limiting, but we postponed giving a justification.  A value $m>0.75$ means that more than 75\% of people infected with the disease die of it; even the Black Death had a mortality fraction less than this.  Worse yet, it is doubtful that the model is appropriate for a disease with such high mortality because of the model assumptions of recruitment and encounter rates that do not depend on the total population.  In reality, a population drastically reduced by disease mortality can be expected to have a much lower birth rate than one not so reduced, and it seems likely that the average number of encounters each of us has each day must be smaller when the population is greatly reduced.  A disease with a mortality fraction higher than 0.75 should be modeled with different assumptions about these processes.

There is no reason why the method presented here should be restricted to epidemiological models.  Any model that has processes with time scales that are significantly different could be approached with the same method.  In a food web model, for example, there could be organisms with a very short life span, such as microorganisms, along with organisms with much longer life spans.  It is common for predator species to have a longer lifespan than their prey, but we should be wary of using the asymptotic approach if the ratio of the lifespans is on the order of 10 rather than 1000, as we generally have in epidemiology models.  But even in these less clear cases, the asymptotic approach can give useful approximations, as long as they are not taken to be quantitatively accurate.

\end{document}